\documentclass{article}
\usepackage{amsmath,amsxtra,amssymb,latexsym,amsfonts,amscd}
\usepackage{multicol,color}
\usepackage{float}
\usepackage{soul}
\usepackage{graphicx}
\usepackage{subfig}
\usepackage{algpseudocode}
\usepackage{mathrsfs}
\usepackage{amssymb}
\usepackage{theorem}
\usepackage{fancyhdr}
\usepackage{tikz}
\usepackage{enumerate}
\usepackage[margin=1.1in]{geometry}
\usepackage{animate}
\usepackage{hyperref}
\usepackage{listings}
\usepackage{xcolor}
\lstdefinestyle{mystyle}{
basicstyle=\small\ttfamily,
breaklines=true,
keywordstyle=\color{blue},
commentstyle=\color{green!40!black},
stringstyle=\color{orange},
numbers=left,
numberstyle=\tiny\color{gray},
frame=single,
columns=fullflexible,
keepspaces=true,
showstringspaces=false,
tabsize=4}
\usetikzlibrary{positioning}
\def\reals{{I\kern-.35em R}}

\usepackage[normalem]{ulem}
\usepackage{bbm}
\def\disp{\displaystyle}

\def\lm{\lambda}
\def\O{\Omega}

\def\({\left(}
\def\){\right)}
\def\[{\left[}
\def\]{\right]}

\def\ox{\bar{x}}

\def\ou{\bar{u}}

\def\gg{\gamma}
\def\dn{\downarrow}

\def\la{\langle}
\def\ra{\rangle}

\def\R{\mathbb{R}}
\def\N{\mathbb{N}}

\def\dn{\downarrow}
\def\O{\Omega}

\def\vph{\varphi}
\def\emp{\emptyset}

\def\lm{\lambda}

\def\gg{\gamma}

\def\al{\alpha}

\def\be{\beta}

\def\ph{\varphi}

\def\N{I\!\!N}
\def\th{\theta}

\def\vt{\vartheta}

\newtheorem{theorem}{Theorem}[section]

\theoremstyle{plain}{\theorembodyfont{\rmfamily}
}
\theoremstyle{plain}{\theorembodyfont{\rmfamily}
}
\theoremstyle{plain}{\theorembodyfont{\rmfamily}
\theoremstyle{plain}{\theorembodyfont{\rmfamily}
}

\theoremstyle{plain}{\theorembodyfont{\rmfamily}
}

\def\eq{\begin{equation}}
\def\eeq{\end{equation}}
\begin{document}
\begin{center}
{\bf OPTIMAL CONTROL OF PERTURBED SWEEPING PROCESSES\\ WITH APPLICATIONS TO GENERAL ROBOTICS MODELS}\\[3ex]
GIOVANNI COLOMBO\footnote{Dipartimento di Matematica ``Tullio Levi-Civita", Universit$\grave{\textrm{a}}$ di Padova, 35121 Padova, Italy (colombo@math.unipd.it), and G.N.A.M.P.A. of Istituto Nazionale di Alta Matematica ``Francesco Severi", Piazzale Aldo Moro 5, 00185 -- Roma, Italy. The research of this author was partly supported by the project funded by the EuropeanUnion – NextGenerationEU under the National Recovery and Resilience Plan (NRRP), Mission 4 Component 2 Investment 1.1 -
Call PRIN 2022 No. 104 of February 2, 2022 of Italian Ministry of
University and Research; Project 2022238YY5 (subject area: PE - Physical
Sciences and Engineering) ``Optimal control problems: analysis,
approximation and applications''.}\quad BORIS S. MORDUKHOVICH\footnote{Department of Mathematics, Wayne State University, Detroit, Michigan 48202, USA (aa1086@wayne.edu). Research of this author was partially supported by the US National Science Foundation under grant DMS-2204519,
by the Australian Research Council under Discovery Project DP-190100555, and by Project~111 of China under grant D21024.} \quad 
DAO NGUYEN\footnote{Department of Mathematics and Statistics, San Diego State University, CA 92182, USA (dnguyen28@sdsu.edu).}\quad 
\quad TRANG NGUYEN\footnote{Department of Mathematics, Wayne State University, Detroit, Michigan 48202, USA (daitrang.nguyen@wayne.edu). Research of this author was partially supported by the US National Science Foundation under grant DMS-2204519.} 
\quad NORMA ORTIZ-ROBINSON\footnote{Department of Mathematics, Grand Valley State University (ortizron@gvsu.edu).}
\end{center}
\small{\sc Abstract.} This paper primarily focuses on the practical applications of optimal control theory for perturbed sweeping processes within the realm of robotics dynamics. By describing these models as controlled sweeping processes with pointwise control and state constraints and by employing necessary optimality conditions for such systems, we formulate optimal control problems suitable to these models and develop numerical algorithms for their solving. Subsequently, we use the Python  Dynamic Optimization library GEKKO to simulate solutions to the posed robotics problems in the case of any fixed number of robots under different initial conditions.\\[1ex]
\normalsize

\section{Introduction and Discussions}\label{applic}
\setcounter{equation}{0}

This paper addresses an important class of {\em differential variational inequalities} (DVIs), which is known as {\em sweeping processes}. This class of discontinuous dynamical systems was introduced and largely investigated in the 1970s by Jean-Jacques Moreau (see, e.g., \cite{M77} and the references therein) via dissipative differential inclusions given in the form
\begin{equation}\label{dynamical system}
\dot{x}(t)\in-N\big(x(t);C(t)\big)\;\mbox{ a.e. }\;t\in[0,T],
\end{equation}
where $C(\cdot)$ is a continuously moving closed and convex set in a Hilbert space, and where $N(\cdot;\Omega)$ for as given set $\Omega$ stands for the classical {\em normal cone} of convex analysis
\begin{equation}\label{nc}
N(x;\O):=\big\{v\;\big|\;\la v,y-x\ra\le 0\;\mbox{ for all }\;y\in \O\big\}\;\textrm{ if }\;x\in\O\;\textrm{ and }\;N(x;\O):=\emp\textrm{ if }\;x\notin\O.
\end{equation}
In view of \eqref{nc}, we see that \eqref{dynamical system} can be equivalently written in the usual DVI form
\begin{equation*}
\la\dot{x}(t),y-x(t) \ra\ge 0\;\mbox{ for all }\;y\in C(t),\;\mbox{ a.e. }\;t\in[0,T].   
\end{equation*}
Over the years, sweeping process theory has been strongly developed and broadly applied to various problems of nonsmooth mechanics, hysteresis, dynamic network, economics, equilibrium systems, etc. We refer the reader to, e.g., \cite{aht,mb,brog,BrS,KP,migorski1,migorski2,MM,St} and the vast bibliographies therein.

One of the most impressive achievements of sweeping process theory is showing that the Cauchy problem $x(0)=x_0$ for \eqref{dynamical system} admits the {\em unique solution} under very natural assumptions; see \cite[Theorem~2]{kunze}. This clearly excludes optimization of the sweeping dynamics described in  \eqref{dynamical system}. Note to this end that the sweeping differential inclusion \eqref{dynamical system}
is significantly different from the model \begin{equation}\label{lip-inc}
\dot{x}(t)\in F\big(x(t)\big)\;\mbox{ a.e. }\;t\in[0,T]
\end{equation}
governed by {\em Lipschitzian} set-valued mappings $F$; in particular, in the case of explicit ODE {{\em control systems} with $F(x):=g(x,U)$, i.e., (considering autonomous systems for simplicity) described by
\begin{equation}\label{ode-control}
\dot{x}(t)=g\big(x(t),u(t)\big),\quad u(t)\in U.   
\end{equation}
Optimal control theory for Lipschitzian differential inclusions \eqref{lip-inc} has been largely developed with deriving comprehensive necessary optimality conditions, etc.; see, e.g., the books \cite{cl,m-book2,v} among numerous of applications. The crucial difference between the developed theory for \eqref{lip-inc} is the set-valued mapping $F(x(t)):=N(x(t);C(t)$ for a.e.\ $t\in[0,T]$ is not just non-Lipschitzian but {\em highly discontinuous}. Furthermore, the uniqueness of solutions to \eqref{dynamical system}, which prevents optimization of such systems, is due to {\em maximal monotonicity} of the unbounded normal cone operators that has never been of interest in \cite{cl,
m-book2,v} and other publications on optimization of Lipschitzian differential inclusions.

Appropriate formulations of optimal control problems for sweeping processes and deriving necessary optimality conditions for sweeping optimal solutions require significant modifications of the original model. Several sweeping control models were proposed and developed in this direction. In \cite{chhm0,chhm}, the authors suggested a {\em control parametrization} of the moving set $C(t)=C(u(t))$ as a controlled halfspace \cite{chhm} or a controlled polyhedron \cite{chhm} with seeking an optimal control that minimises a cost functional of the Mayer or Bolza type. Another approach to sweeping optimal control was suggested and developed in \cite{bk}, where control functions entered the {\em adjacent ODE system}; see also \cite{ao} for further results in this vein. 

However, the main attention in the most recent developments on sweeping optimal control has been paid to control systems of the type
\begin{equation}\label{Problem}
\left\{\begin{matrix}
\dot{x}(t)\in-N\big(x(t);C\big)+g\big(x(t),u(t)\big)\;\textrm{ a.e. }\;t\in[0,T],\;x(0)=x_0\in C\subset\R^n,\\
u(t)\in U\subset\R^d\;\textrm{ a.e. }\;t\in[0,T],
\end{matrix}\right.
\end{equation}
where {\em control actions} $u(t)$ appear in the {\em dynamic perturbations} $g(x,u)$. Sweeping control systems of type \eqref{Problem} are direct extensions of the standard controlled ODE dynamics in \eqref{ode-control} by adding the normal cone term $N(x(t);C)$, which generates completely new phenomena when $\emp\ne C\ne\R^n$. Observe the following:

$\bullet$ The {\em feasible velocity set}
\begin{equation}\label{velocity}
F(x):=-N(x;C)+g(x,U),\;\mbox{ where }\;g(x,U):=\big\{v\in\R^n\;\big|\;\exists\,u\in U\;\mbox{ with }\;v=g(x,u)\big\},
\end{equation}
is {\em unbounded} and {\em discontinuous} for all $x\in C$.

$\bullet$ There are intrinsic {\em pointwise state constraints}
\begin{equation}\label{state}
x(t)\in C\;\mbox{ for all }\;t\in[0,T].
\end{equation}
Indeed, \eqref{state} immediately follows from \eqref{Problem} due to the fulfillment of $N(x(t);C)\ne\emp$ on $[0,T]$ therein. 

To the best of our knowledge, a sweeping optimal control model with controlled perturbations of type \eqref{Problem} has been formulated and studied for the first time in \cite{CM16}, in the case where $C$ is a convex polyhedron, with establishing there necessary conditions for optimal solutions by using the method of {\em discrete approximations} developed in \cite{m95,m-book2} for Lipschitzian differential inclusions and then significantly extended to controlled sweeping processes in many publications. Over the recent years, various methods to derive necessary optimality conditions of different types for sweeping control problems with controlled perturbations have been proposed and implemented under diverse assumptions on the sweeping set $C$ and the perturbation mapping $g$; see \cite{CCMN21,CKMNP21a,b1,cg21,CMN20,cmnn23b,pfs23,her-pall,zeidan} among other publications. The obtained optimality conditions were applied to solving a variety of practical models: crowd motions and traffic equilibria \cite{CM17,b1,CMN19}, locomotion of a soft-robotic crawler \cite{cg21}, mobile robot dynamics \cite{CMN19}, marine surface vehicle modeling \cite{CKMNP21a,mnn23a}, nanoparticle modeling \cite{mnn23a}, etc.
 
This paper continue the lines of new applications of sweeping control theory to practical modeling. Our main attention is paid to a {\em general robotic model} with obstacles, which is formulated here as an optimal control problem for the controlled sweeping dynamics \eqref{Problem} with pointwise state and control constraints. A much simplified version of this model was considered in \cite{CMN19} based on necessary optimality conditions derived in \cite{CMN20}. The results obtained in \cite{CMN19,CMN20} allowed us to determine an optimal control only for the model therein with only one obstacle. Now, by using newly developed necessary optimality conditions taken from \cite{cmnn23b}, we are able to solve the controlled robotic model for an arbitrarily fixed number of obstacles, design a numerical algorithm to find optimal solutions, and constructively implement it using the Python Dynamic Optimization library via GEKKO Optimization Suite. 

The rest of the paper is organized as follows. In Section~\ref{sec:sweep}, we formulate an optimal control problem for the control sweeping dynamic \eqref{Problem} and present the necessary optimality conditions for solutions to this problem taken from \cite{cmnn23b} and applied below to the robotic model. The controlled robotic model of our study is described and discussed in Section~\ref{Prob-Form}. The subsequent Section~\ref{NOC} applies the obtained necessary optimality conditions to the controlled sweeping robotic model. Based on these optimality condition, we design in Section~\ref{Num} the algorithm to solve the controlled robotic model and implement it numerically by using GEKKO in various model settings and scenario to find optimal solutions.

\section{Optimization of Controlled Sweeping Processes}\label{sec:sweep}

In this section we present, for completeness and reader's convenience, the sweeping optimal control problem and necessary conditions for its optimal solutions taken from \cite{cmnn23b}. In contrast to \cite{cmnn23b}, we consider here a control problem on a fixed time interval $[0,T]$, which is used on our robotics application. Moreover, we do not need to impose additional endpoint constraints $x(T)\in\O$ as in \cite{cmnn23b}.

The precise formulation of the sweeping optimal control problem $(P)$ is as follows: minimize the Mayer-type cost functional
\begin{equation}\label{cost1}
J[x,u]:=\varphi\big(x(T)\big)
\end{equation}
over control functions $u(\cdot)$ and the corresponding trajectories $x(\cdot)$ satisfying the controlled sweeping dynamics in \eqref{Problem} with pointwise control constraints, where the set $C$ is a {\em convex  polyhedron} given by
\begin{equation}\label{C}
\begin{array}{ll}
C:=\bigcap_{j=1}^{s}C^j\textrm{ with }C^j:=\big\{x\in\R^n\;\big|\;\la x^j_*,x\ra\le c_j\big\},
\\
\|x_\ast^j\|=1$, $j=1,\ldots,s.
\end{array}
\end{equation}
As follows from \eqref{state} and \eqref{C}, we automatically have the pointwise {\em state constraints}
\begin{equation}\label{stateC}
\la x^j_*,x(t)\ra\le c_j\;\mbox{ for all }\;t\in[0,T]\;\mbox{ and }\;j=1,\ldots,s.
\end{equation}
The pair $(x(\cdot),u(\cdot))$ satisfying \eqref{Problem} and \eqref{stateC} is a {\em feasible solution} to $(P)$ if $u(\cdot)\in L^2([0,T];\R^d)$ and $x(\cdot)\in W^{1,2}([0,T];\R^n)$. We say that a feasible solution $(\ox(\cdot),\ou(\cdot))$ is {\em $W^{1,2}\times L^2$-local minimizer} in $(P)$ if $J[\ox,\ou]\le J[x,u]$ for all feasible solutions $(x(\cdot),u(\cdot))$ satisfying the localization condition
\begin{equation*}
\int_0^T\Big(\|\dot{x}(t)-\dot\ox(t)\|^2+\|u(t)-\ou(t)\|^2\Big)dt\le\varepsilon.
\end{equation*}
From now on, we call the above pair $(\ox(\cdot),\ou(\cdot))$ simply a ``local minimizer" in $(P)$ and observe that this notion of local minima is specific for variational problems of type $(P)$. It clearly holds for {\em strong minimizers} with the replacement of $L^2$ by the space ${\cal C}$ and, of course, for global optimal solutions to $(P)$.   

Next we formulate the basic assumptions needed for fulfillment of the necessary optimality conditions presented bellow for the $W^{1,2}\times L^2$-local minimizer $(\ox(\cdot),\ou(\cdot))$ in $(P)$ under consideration. Observe that the main results of \cite{cmnn23b} are established under essentially less restrictive assumptions while those formulated below are met in the robotics model of our interest in this paper.   

{\bf(H1)} The control region $U\ne\emp$ is a closed, bounded, and convex subset of $\R^d$.

{\bf(H2)} The perturbation mapping $g\colon\R^n\times U\to\R^n$ is continuously differentiable around $\ox(t),\ou(t))$ on $[0,T]$. Furthermore, $g$ satisfies the the sublinear growth condition
\begin{equation*}
\|g(x,u)\|\le\be\big(1+\|x\|\big)\;\mbox{ for all }\;u\in U\;
\mbox{ with some }\;\be>0.
\end{equation*}

{\bf(H3)} The vector field $g(x,U)$ from \eqref{velocity} is convex for all $x$ around $\ox(t)$ and all $t\in[0,T]$.

{\bf(H4)} The linear independence constraint qualification (LICQ) is satisfied along $\ox(t)$ on $[0,T]$, i.e.,
\begin{equation*}
\Big[\sum_{j\in I(\ox)}\al_jx^j_*=0,\;\al_j\in\R\Big]\Longrightarrow\big[\al_j=0\;
\text{ for all }\;j\in I(\ox)\big],
\end{equation*}
where $I(x)$ stands for the collections of active indices of polyhedron \eqref{C} at $x$ defined by
\begin{equation}\label{I}
I(x):=\big\{j\in\{1,\ldots,s\}\;\big|\;\la x^j_*,x\ra=c_j\big\}. 
\end{equation}

{\bf(H5)} The cost function $\ph$ is continuously differentiable around $\ox(T)$. 

Deduce further from the Motzkin's theorem of the alternative that the normal cone to the polyhedral set $C$ from \eqref{C} is represented in the form
\begin{equation}\label{alter}
N(x;C)=\Big\{\sum_{j\in I(x)}\eta^j x^j_*\;\Big|\;\eta^j\ge 0\Big\}.
\end{equation}

Using \eqref{alter}, it follows that any trajectory $x(\cdot)$ of \eqref{Problem} is represented as  
$$
-\dot{x}(t)=\sum_{j=1}^s\eta^j(t)x^j_*-g\big(x(t),u(t)\big)\;\textrm{ for all }\;t\in [0,T),
$$
where $\eta^j\in L^2([0,T];\mathbb{R}^+)$ and $\eta^j(t) =0$ for all $t$ such that $j\notin I(t,x(t))$. We say that the {\em normal cone to $C$ is active} along $x(\cdot)$ on the set $E\subset [0,T]$ if for a.e. $t\in E$ and all $j\in I(x)$ it holds
$\eta^j(t)>0$. Denoting by $E_0$ the largest subset $E$ of $[0,T]$ where the normal cone to $C$ is active along $\bar{x}$ on $E$, we say that it is active along $x(\cdot)$ provided that $E_0=[0,T]$.

Now we are ready to formulate the necessary optimality conditions for problem $(P)$, which specify the results of \cite[Theorem~5.2]{cmnn23b} for problem $(P)$ under consideration.

\begin{theorem}\label{Thm6.1*}
Let $(\ox(\cdot),\ou(\cdot))$ be a local minimizer to problem $(P)$ under the imposed assumptions {\rm(H1)--(H5)}. Then there exist a multiplier $\lm\ge 0$, a nonnegative vector measure 
$\gg_>=(\gg_>^1,\ldots,\gg_>^s)\in C^*([0,T];\R^s)$, a signed vector measure $\gg_0=(\gg_0^1,\ldots,\gg_0^s)\in C^*([0,T];\R^s)$, as well as adjoint arcs 
$p(\cdot)\in W^{1,2}([0,T];\R^n)$ and $q(\cdot)\in BV([0,T];\R^n)$ such that the following conditions are fulfilled:\\[1ex]
$\bullet$ The {\sc primal arc representation}
\begin{equation*}
-\dot{\ox}(t)=\sum_{j=1}^s\eta^j(t)x^j_*-g\big(\ox(t),\ou(t)\big)\;\textrm{ for a.e. }\;t\in [0,T),
\end{equation*}
where the functions $\eta^j(\cdot)\in L^2([0,T);\R_+)$, $j=1,\ldots,s$, are uniquely determined for a.e.\ $t\in[0,T)$ by the above representation.\\[1ex]
$\bullet$ The {\sc adjoint dynamical system}
\begin{equation*}
\dot{p}(t)=-\nabla_x g\big(\ox(t),\ou(t)\big)^*q(t)\;\textrm{ for a.e. }\;t\in[0,T],
\end{equation*}
where the right continuous representative of $q(\cdot)$, with the same notation, satisfies
\begin{equation*}
q(t)=p(t)-\int_{(t,T]} \sum_{j=1}^sd\gg^j(\tau)x^j_*
\end{equation*}
for all $t\in[0,T]$ except at most a countable subset, and moreover $p(T)=q(T)$ with $\gamma=\gamma_>+\gamma_0$.\\[1ex]
$\bullet$ The maximization condition
\begin{equation*}
\big\la\psi(t),\ou(t)\big\ra=\max_{u\in U}\big\la\psi(t),u\big\ra\;\textrm{ for a.e. }\;t\in[0,T],
\end{equation*}
$\bullet$ The {\sc dynamic complementary slackness conditions}
\begin{equation*}
\big\la x^j_*,\ox(t)\big\ra<c_j\Longrightarrow\eta^j(t)=0\;\mbox{and }\;\eta^j(t)>0\Longrightarrow\;\big\la x^j_*,q(t)\big\ra=0
\end{equation*}
for a.e. $t\in[0,T]$ and all indices $j=1,\ldots,s$.\\[1ex]
$\bullet$ The {\sc transversality condition}: there exist numbers $\eta^j(T)\ge 0$ for $j\in I(\ox(T))$ such that
\begin{equation*}
-p(T)=\sum_{j\in I(\ox(T))}\eta^j(T)x^j_*+\lm\nabla\vph\big(\ox(T)\big)\;\textrm{ and }\;\sum_{j\in I(\ox(T))}\eta^j(T)x^j_*\in
N\big(\ox(T);C\big),
\end{equation*}
where the collection of active indices of the polyhedron \eqref{C} is taken from \eqref{I}.\\[1ex] 
$\bullet$ The {\sc endpoint complementary slackness conditions}
\begin{equation*}
\big\la x^j_*,\ox(T)\big\ra<c_j\Longrightarrow\eta^j(T)=0\quad\mbox{for all }\;j\in I\big(\ox(T)\big).
\end{equation*}
$\bullet$ The {\sc measure nonatomicity condition:} If $t\in[0,T)$ and $\la x^j_*,\ox(t)\ra<c_j$ for all $j=1,\ldots,s$, then there
exists a neighborhood $V_t$ of $t$ in $[0,T)$ such that $\gg(V)=0$ for all the Borel subsets $V$ of $V_t$.\\[1ex]
$\bullet$ The  {\sc general nontriviality condition}\\[1ex]
\begin{equation*}
(\lm,p,\|\gamma_0\|_{TV},\|\gg_>\|_{TV})\ne 0
\end{equation*}
accompanied by the {\sc support condition}\\[1ex]
\begin{equation*}
\mathrm{supp}(\gamma_>)\cap \mathrm{int}(E_0)=\emp,
\end{equation*}
which holds provided that the normal cone is active on a set with nonempty interior.\\[1ex]
$\bullet$ The {\sc enhanced nontriviality condition}
\begin{equation*}
(\lm,p(T))\ne 0 
\end{equation*}
provided that $\la x^j_*,\ox(t)\ra<c_j$ for all $t\in[0,T]$ and all indices $j=1,\ldots,s$.
\end{theorem}

The proof of Theorem~\ref{Thm6.1*}, as induced by that of \cite[Theorem~5.2]{cmnn23b}, constitutes a special case where the time is fixed. It relies on the {\em method of discrete approximations} and the machinery of {\em variational analysis} developed in \cite{m-book1,m24}. Note that,
although problem $(P)$ is formulated only in terms of smooth functions and convex sets, the given proof employs tools of {\em second-order generalized differentiation} involving the coderivative calculation for the feasible velocity mapping \eqref{velocity} whose graph is always {\em nonconvex}; see \cite{cmnn23b,m24} for more details and references.

\section{Controlled Robotic Model}\label{Prob-Form}
\setcounter{equation}{0}

This section is dedicated to recalling the formulation of an optimal control problem for the robotics model, featuring obstacles, whose dynamics are determined as a sweeping process. Generally, the model considers $n\ge 2$ robots represented as safety disks with the same radius $R$ located on a plane. Each robot aims to reach its destination via the shortest possible path during the designated time interval $[0,T]$ while avoiding collisions (but potentially making contact) with other obstacles, such as walls or other robots. 

First we describe the trajectory $x^i(t)$ of the $i$-robotic by
\begin{eqnarray*}
x^i(t)=\big(\|x^i(t)\|\cos\th_i(t),\|x^i(t)\|\sin\th_i(t)\big)\;\mbox{ for }\;i=1,\ldots,n
\end{eqnarray*}
with the configuration vector $x=(x^1,\ldots,x^n)\in\R^{2n}$, where $x^i\in\R^2$ denotes the center of the $i$-th robotic disk with coordinates representation $(\|x^i\|\cos\th_i,\|x^i\|\sin\th_i)$ and $\th_i$ stands for the smallest positive angle in standard position formed by the positive $x$-axis and $Ox^i$, with the target is the origin. 

To avoid collisions during the movement of all robots, the distance between the robots and obstacles must be maintained. Hence the set of {\em admissible configurations} is formulated by imposing the {\em nonoverlapping conditions} $\|x^i-x^j\|\ge 2R$ as follows:
\begin{eqnarray}\label{ac}
A:=\big\{x=\(x^1,\ldots,x^n\)\in\mathbb{R}^{2n}\big|\;D_{ij}(x)\ge 0\;\mbox{ for all }\;i,j\in\{1,\ldots,n\}\big\},
\end{eqnarray}
where $D_{ij}(x)=\|x^{i}-x^j\|-2R$ is the distance between the inelastic disks $i$ and $j$. Let $\nabla D_{ij}(x)$ be the gradient of the distance function $D_{ij}(x)$ at $x\ne 0$, the set of {\em admissible velocities} is defined by
\begin{eqnarray*}
V_h(x):=\big\{v\in\R^{2n}\big|\;D_{ij}(x)+h\nabla D_{ij}(x)v\ge 0\;\mbox{ for all }\;i,j\in\{1,\ldots,n\},\;i<j\big\},\quad x\in\R^{2n}.
\end{eqnarray*}
Suppose that at time $t_k\in[0,T]$ the  admissible configuration is $x_k:=x(t_k)\in A$. Then after the period of time $h>0$, the next configuration is $x_{k+1}=x(t_k+h)$. Taking into account the first-order Taylor expansion around $x_k\ne 0$ gives us 
\begin{eqnarray}\label{Dij}
D_{ij}\big(x(t_k+h)\big)=D_{ij}\big(x(t_k)\big)+h\nabla D_{ij}\big(x(t_k)\big)\dot{x}(t_k)+o(h)\;\mbox{ for small }\;h>0.
\end{eqnarray}
Considering the admissible velocity $\dot{x}(t_k)\in V_h(x_k)$ and skipping the term $o(h)$ for small $h$ yields
\begin{eqnarray*}
D_{ij}(x_k)+h\big\la\nabla D_{ij}(x_k),\dot{x}(t_k)\big\ra\ge 0,
\end{eqnarray*}
and therefore it follows from \eqref{Dij} that 
$$
D_{ij}(x(t_k+h))\ge 0, \mbox{ i.e.},\, x(t_k+h)\in A.
$$

In the absence of other obstacles, every robot aims to reach its target by the {\em shortest path}; so the robots tend to keep their desired spontaneous velocities till reaching the target. However, when the robot in consideration touches the obstacles in the sense that $\|x^{i}(t)-x^j(t)\|=2R$, its velocity needs adjustment to ensure the distance remains at least $2R$ through some specific {\em control actions} applied to the velocity term. This situation is represented through the following modeling approach:
\begin{eqnarray}\label{g}
g\big(x(t),u(t)\big)=\big(s_1\|u^1(t)\|\cos\th_1(t),s_1\|u^1(t)\|\sin\th_1(t),\ldots,s_n\|u^n(t)\|\cos\th_n(t),s_n\|u^n(t)\|\sin\th_n(t)\big)
\end{eqnarray}
where the control constraints are defined by
\begin{eqnarray}\label{u}
u(t)=\big(u^1(t),\ldots,u^n(t)\big)\in U\subset\R^n\;\mbox{ for a.e. }t\in[0,T]
\end{eqnarray}
with the closed and convex control region $U$ specified below.

To prevent the overlap between the robot and obstacles, let $x_k\in A$ at the time $t_k$, and use the mapping $g\colon\R^{2n}\times\R^n\to\R^{2n}$ from \eqref{g} with a given feasible control $u_k:=u(t_k)$ from \eqref{u}. The next configuration $x_{k+1}$ is computed by
\begin{eqnarray}\label{dsg}
x_{k+1}=x_k+h V_{k+1},\;V_{k+1}\in\R^{2n},
\end{eqnarray}
addressing the following {\em convex optimization problem}:
\begin{eqnarray}\label{P_i}
\mbox{minimize }\;\|g(x_k,u_k)-V\|^2\;\mbox{ subject to }\;V\in V_{h}(x_k),
\end{eqnarray}
where the control $u_k\in U$ is incorporated into the desired velocity term to modify the actual velocities of the robots, guaranteeing they maintain nonoverlapping trajectories. Moreover, the algorithm described in equations \eqref{dsg} and \eqref{P_i} implies that  $V_{k+1}$ is chosen as the (unique) element from the admissible velocities set that is nearest to the desired velocity $g(x_k,u_k)$, preventing thereby any overlap among the robots. Select any natural number $m$ and portion the interval $[0,T]$ into the $2^m$ equal subintervals of length $h_m:=T/2^m\dn 0$ as $m\to\infty$. Considering the discrete time $t_{km}:=kh_m$, denote $I_{km}:=[t_{km},t_{(k+1)m})$ for $k=0,\ldots,2^m-1$ and $I_{2^m m}:=\{T\}$. Hence, by referring subsequently to equations  \eqref{dsg} and \eqref{P_i}, the algorithm is formulated as follows:
$$x_{0m}\in A$$
\begin{eqnarray}\label{alg}
x_{(k+1)m}:=x_{km}+h_m V_{(k+1)m},\;k=0,\ldots,2^m-1,
\end{eqnarray}
where the admissible velocity in the next iteration is defined by
\begin{eqnarray}\label{V}
V_{(k+1)m}:=\disp\Pi\big(g(x_{km},u_{km});V_{h_m}(x_{km})\big),\, k=0,\ldots,2^m-1.
\end{eqnarray}
Taking into account the construction of $x_{km}$ for $0\le k\le 2^m-1$ and $m\in\N$, define a sequence of piecewise linear mappings $x_{2^m}\colon[0,T]\to\R^{2n}$, $m\in\N$, which passes through those points by:
\begin{eqnarray}\label{x_k}
x_{2^m}(t):=x_{km}+(t-t_{km})V_{(k+1)m}\;\mbox{ for all }\;t\in I_{km},\quad k=0,\ldots,2^m-1.
\end{eqnarray}
The following relationships hold whenever $m\in\N$:
\begin{eqnarray}\label{x_2m}
x_{2^m}(t_{km})=x_{km}=\underset{t\to t_{km}}{\lim}x_{km}(t)\;\mbox{ and }\;\dot{x}_{2^m}(t):=V_{(k+1)m}\;\mbox{ for all }\;t\in(t_{km},t_{(k+1)m}).
\end{eqnarray}
Through discussions in \cite{mb1}, it is observed that the solutions to \eqref{x_k} in the {\em uncontrolled} setting of \eqref{V} with $g=g(x)$ uniformly converge on $[0,T]$ to a trajectory of a specific perturbed sweeping process.  The {\em controlled} model under investigation in this scenario presents significantly more involved. To advance and utilize the insights from \cite[Theorem~5.2]{cmnn23b}, the following set is considered for all $x\in\R^{2n}$ 
\begin{eqnarray}\label{K}
K(x):=\big\{y\in\R^{2n}\big|\;D_{ij}(x)+\nabla D_{ij}(x)(y-x)\ge 0\;\mbox{ whenever }\;i<j\big\},
\end{eqnarray}
enabling us to express the algorithm in \eqref{V}, \eqref{x_k} as
\begin{eqnarray*}
x_{(k+1)m}=\Pi\Big(x_{km}+h_m g(x_{km},u_{km});K(x_{km})\Big)\;\mbox{ for }\;k=0,\ldots,2^m-1.
\end{eqnarray*}
It then can be expressed in an equivalent form as
\begin{eqnarray*}
x_{2^m}\big(\vt_{2^m}(t)\big)=\Pi\Big(x_{2^m}(\tau_{2^m}(t))+h_m g\big(x_{2^m}(\tau_{2^m}(t)),u_{2^m}(\tau_{2^m}(t)\big);K(x_{2^m}(\tau_{2^m}(t))\Big)\;\mbox{ for all }\;t\in[0,T],
\end{eqnarray*}
where the functions $\tau_{2^m}(\cdot)$ and $\vt_{2^m}(\cdot)$ are extended to the entire interval $[0,T]$ by $\tau_{2^m}(t):=t_{km}$ and $\vt_{2^m}(t):=t_{(k+1)m}$ for all $t\in I_{km}$. In addition, considering the construction of the convex set $K(x)$ as described in \eqref{K} and the definition of the normal cone \eqref{nc}, along with the relationships outlined in \eqref{x_2m}, we get the following controlled sweeping dynamics:
\begin{eqnarray}\label{dotx}
\dot x_{2^m}(t)\in-N\big(x_{2^m}(\vt_{2^m}(t));K(x_{2^m}(\tau_{2^m}(t)))\big)+g\big(x_{2^m}(\tau_{2^m}(t)),u_{2^m}(\tau_{2^m}(t))\big)\;\mbox{ a.e. }\;t\in[0,T]
\end{eqnarray}
with $x_{2^m}(0)=x_0\in K(x_0)=A$ and $x_{2^m}(\vt_{2^m}(t))\in K(x_{2^m}(\tau_{2^m}(t)))$ on $[0,T]$. The inclusion in \eqref{dotx} can be formalized as a controlled perturbed sweeping process \eqref{Problem} through the convex polyhedron 
\begin{eqnarray}\label{2.13}
C:=\bigcap_{j=1}^{n-1}\big\{x\in\R^{2n}\big|\;\la x^j_*,x\ra\le c_j\big\}
\end{eqnarray}
with $c_j:=-2R$ and with the $n-1$ vertices of the polyhedron given by
\begin{eqnarray}\label{e}
x^j_*:=e_{j1}+e_{j2}-e_{(j+1)1}-e_{(j+1)2}\in\R^{2n},\quad j=1,\ldots,n-1,
\end{eqnarray}
where $e_{ji}$ for $j=1,\ldots,n$ and $i=1,2$ are the vectors in $\R^{2n}$ of the form
\begin{eqnarray*}
e:=\big(e_{11},e_{12},e_{21},e_{22},\ldots,e_{n1},e_{n2}\big)\in\R^{2n}
\end{eqnarray*}
with 1 at only one position of $e_{ji}$ and $0$ at all the other positions.

As a result, we proceed to formulate {\em sweeping optimal control problem} denoted as (P) from Section~2, characterized as the continuous-time counterpart of the discrete algorithm employed in the robotics model. Consider the cost functional
\begin{equation}\label{t:102*}
\mbox{minimize }\;J[x,u]:=\disp\frac{1}{2}\big\|x(T)\big\|^2,
\end{equation}
with the model objective of  {\em minimizing the distance} of robots from the admissible configuration set to the target. We describe the continuous-time dynamics by the controlled sweeping process
\begin{eqnarray}\label{t:101*}
\left\{\begin{array}{lcl}
-\dot{x}(t)\in N\big(x(t);C\big)+g\big(x(t),u(t)\big)\;\mbox{ for a.e. }\;t\in[0,T],\\
x(0)=x_0\in C,\;u(t)\in U\;\mbox{ a.e. on }[0,T],
\end{array}\right.
\end{eqnarray}
where the constant set $C$ is taken from \eqref{2.13}, the control constraints reduce to \eqref{g}. The dynamic requirement for nonoverlapping, expressed as $\|x^i(t)-x^j(t)\|\ge 2R$, is equivalent to the pointwise state constraints given by
\begin{eqnarray}\label{t:100*}
x(t)\in C\Longleftrightarrow\la x^j_*,x(t)\ra\le c_j\;\mbox{ for all }\;t\in[0,T]\;\mbox{ and }\;j=1,\ldots,n-1,
\end{eqnarray}
which follows from \eqref{x_2m} s a result of constructing $C$ and the normal cone definition. Additionally, considering that the result in the limiting process from Theorem~\ref{Thm6.1*} will be employed when the discrete step $h$ in \eqref{dsg} diminishes, then for convenience, we choose an equivalent norm $\|(x^{j1},x^{j2})\|:=|x^{j1}|+|x^{j2}|$ for each component $x^j\in\R^2$ of $x\in\R^{2n}$.

\section{Necessary Optimality Conditions for the Robotics Model}\label{NOC}\setcounter{equation}{0}

In this section, we obtain the necessary optimality conditions for problem defined in \eqref{t:102*}--\eqref{t:100*} with the robotics model data by applying Theorem~\ref{Thm6.1*}. However, the absence of the varying time $T$ can be easily incorporated into the proof, and so we omit it here while referring the reader to \cite{CMN19}.

\begin{theorem}{\bf(necessary optimality conditions for the sweeping controlled robotics  model).}\label{Thm6.2*a} 
Let $(\ox(\cdot),\ou(\cdot))$ be a $W^{1,2}\times L^2$-local minimizer for problem $(P)$, and let all the assumptions of Theorem~{\rm\ref{Thm6.1*}} be fulfilled.
Then there exist a multiplier $\lm\ge 0$, a measure $\gg\in C^*([0,T];\R^{2n})$ as well as adjoint arcs $p(\cdot)\in W^{1,2}([0,T];\R^{n})$ and $q(\cdot)\in BV([0,T];\R^{n})$ satisfying to the conditions:
\begin{itemize}
\item[\bf(1)] $-\dot{\ox}(t)=\disp{\sum^{n-1}_{j=1}}\eta^j(t)x_*^j+\big(g(\ox^1(t),\ou^1(t))),\ldots,g(\ox^n(t),\ou^n(t)\big)$ for
a.e.\ $t\in[0,T)$,\\ where $\eta^j(\cdot)\in L^2([0,T];\R_+)$ are uniquely defined by this representation and well defined at $t=T$.
\item[\bf(2)] $\|\ox^{j}(t)-\ox^{j+1}(t)\|>2R\Longrightarrow\eta^j(t)=0$ for all $j=1,\ldots,n-1$ and a.e.\ $t\in[0,T]$ including $t=T$.
\item[\bf(3)] $\eta^j(t)>0\Longrightarrow\la x^j_*,q(t)\ra=0$ for all $j=1,\ldots,n-1$ and a.e.\ $t\in[0,T]$.
\item[\bf(4)] $\dot{p}(t)=-\nabla_x g\big(\ox(t),\ou(t)\big)^{\ast} q(t)$ for a.e. $t\in[0,T]$.
\item[\bf(5)] $q(t)=p(t)-\gamma([t,T])$ for all $t\in[0,T]$ except at most a countable subset.
\item[\bf(6)] $\big\la\psi(t),\ou(t)\big\ra=\max_{u\in U}\big\la\psi(t),u\big\ra\;\mbox{ for a.e. }\;t\in[0,T]$.
\item[\bf(7)] $-p(T)=\lm\ox(T)+\sum_{j\in I(\ox(T))}\eta^j(T)x^j_*$ via the set of active polyhedron indices $I(\ox(T))$ at $\ox(T)$.
\item[\bf(8)] $\sum_{j\in I(\ox(T))}\eta^j(T)x^j_*\in N\big(\ox(T);C)$.
\item[\bf(9)] $(\lm,p,\|\gamma_0\|_{TV},\|\gg_>\|_{TV})\ne 0$  with the support conditions
\begin{equation*}
\left\{\begin{array}{ll} 
\mathrm{supp}(\gg_0)\cup\;\mathrm{supp}(\gg_>)\subset
\{t\;|\:\|\ox^{j}(t)-\ox^{j+1}(t)\|=2R\},\\
\mathrm{supp}(\gg_>)\cap \mathrm{int}(E_0)=\emptyset,\;\textrm{ provided that }\;\mathrm{int}(E_0)=\mathrm{int}(\{t\;|\;\eta^j(t)>0, j\in I(x)\})\;\textrm{ is nonempty}. 
\end{array}
\right.
\end{equation*}
\end{itemize}
\end{theorem}

Subsequently, the above theorem allows us to deduce several conclusions for the robotics model, considering the perturbation mapping $g$  as described in equation \eqref{g}.

$\bullet$ It can be seen from the model description that if at he contacting time $t_1\in[0,T]$ we have $\|\ox^{i}(t_1)-\ox^1(t_1)\|=2R$ for some $i=2,\ldots,n$, then the considered robot tends to adjust its velocity to ensure a distance of at least $2R$
between itself and the obstacle in contact with. According to the model requirements, the robot will keep a constant velocity after time $t=t_1$ until either reaching other obstacles ahead or stopping at the ending time $t=T$. Moreover, if the robot comes into contact with other obstacles, it pushes them to move toward the end of the straight circular capillaries in the same direction as before $t=t_1$. Using \eqref{g}, we reformulate the differential relation in {\bf(1)} for a.e.\ $t\in[0,T]$ as follows:
\begin{eqnarray}\label{t:103*}
\left\{\begin{array}{ll}
-\big(\dot{\ox}^{11}(t),\dot{\ox}^{12}(t)\big)=\eta^1(t)(1,1)-\big(s_1\ou^1(t)\cos\th_1(t),s_1\ou^1(t)\sin\th_1(t)\big),\\
-\big(\dot{\ox}^{i1}(t),\dot{\ox}^{i2}(t)\big)=\eta^{i-1}(t)(-1,-1)+\eta^{i}(t)(1,1)-\big(s_i\ou^i(t)\cos\th_i(t),s_i\ou^i(t)\sin\th_i(t)\big)\\
\mbox{ whenever }\;i=2,\ldots,n-1,\mbox{ and}\\
-(\dot{\ox}^{n1}(t),\dot{\ox}^{n2}(t)\big)=\eta^{n-1}(t)(-1,-1)-\big(s_n\ou^n(t)\cos\th_n(t),s_n\ou^n(t)
\sin\th_n(t)\big).
\end{array}\right.
\end{eqnarray}

$\bullet$ If the robot under consideration (referred to as robot 1) does not make contact with the first obstacle (robot 2) in the sense that $\|\ox^2(t)-\ox^1(t)\|>2R$ for all $t\in[0, T]$, then we deduce from {\bf(2)} of Theorem~\ref{Thm6.2*a} that $\eta^1(t)=0$ for a.e. $t\in[0,T]$.
Substituting $\eta^1(t)=0$ into \eqref{t:103*} yields 
\begin{eqnarray*}
-\big(\dot{\ox}^{11}(t),\dot{\ox}^{12}(t)\big)=-\big(s_1\ou^1(t)\cos\th_1(t),s_1\ou^1(t)\sin\th_1(t)\big)\;\mbox{ a.e. on }\;[0,T],
\end{eqnarray*}
which means that the actual velocity and the spontaneous velocity of the robot agree for a.e.\ $t\in[0,T]$. Similarly we conclude that the condition $\|\ox^n(t)-\ox^{n-1}(t)\|>2R$ on $[0,T]$ yields $-\dot{\ox}^n(t)=-g(\ox^n(t),\ou^n(t))$ for a.e.\ $t\in[0,T]$, and then continue in this way with robot $i$.

To proceed further, assume that $\lm=1$ (otherwise we do not have enough information to efficiently employ Theorem~\ref{Thm6.2*a}). Additionally, for simplicity in handling the forthcoming examples, suppose that $\ou^i(\cdot)$ are constant on $[0,T]$ for all $i=1,\ldots,n$. Applying the Newton-Leibniz formula in \eqref{t:103*}, we obtain the trajectory representations
\begin{eqnarray}\label{nlf}
\left\{\begin{array}{ll}
\(\ox^{11}(t),\ox^{12}(t)\)=\(x^{11}_0,x^{12}_0\)-\disp\int^t_0\eta^1(\tau)\(1,1\)d\tau+t\(s_1\ou^1\cos\th_1,s_1\ou^1\sin\th_1\),\\
\(\ox^{i1}(t),\ox^{i2}(t)\)=\(x^{i1}_0,x^{i2}_0\)+\disp\int^t_0\eta^{i-1}(\tau)\(1,1\)d\tau-\disp\int^t_0\eta^{i}(\tau)\(1,1\)d\tau\\
\quad\quad\quad\quad\quad\quad\quad+t\(s_i\ou^i\cos\th_i,s_i\ou^i\sin\th_i\)\;\mbox{ whenever }\;i=2,\ldots,n-1,\\
\(\ox^{n1}(t),\ox^{n2}(t)\)=\(x^{n1}_0,x^{n2}_0\)+\disp\int^t_0\eta^{n-1}(\tau)\(1,1\)d\tau+t\(s_n\ou^n\cos\th_n,s_n\ou^n\sin\th_n\)
\end{array}\right.
\end{eqnarray}
for all $t\in[0,T]$, where $x_0:=(x^{11}_0,x^{12}_0\ldots,x^{n1}_0,x^{n2}_0)\in C$ stands for the starting point in \eqref{t:101*}.

As follows from condition $\textbf{(2)}$, we can assume that $\eta(\cdot)$ is piecewise constant on $[0,T]$ and satisfies
$$\eta(t)=0\;\text{ for a.e. }\;t\in[0,t_i)\;\text{ and }\;\eta(t)=\eta^i(t)\;\text{ for a.e. }\;t\in[t_i,T],$$
where $t_i$ is the contacting time between robot $i$ and robot $i+1$. Using condition \ref{t:103*}, the velocities of robots $i$ and $i+1$ after the contacting time can be expressed as
\begin{eqnarray}
\left\{\begin{array}{ll}
\dot{\ox}^{i}(t)=-\eta^i(t)(1,1)+\big(s_i\ou^i(t)\cos\th_i(t),s_i\ou^i(t)\sin\th_i(t)\big),\\
\dot{\ox}^{i+1}(t)=-\eta^i(t)(-1,-1)+\big(s_{i+1}\ou^{i+1}(t)\cos\th_{i+1}(t),s_{i+1}\ou^{i+1}(t)
\sin\th_{i+1}(t)\big)
\end{array}\right.
\end{eqnarray}
for $i=1,\ldots, n-1$. Hence we obtain the following representation 
$$\eta^i(t)=\begin{cases}\frac{1}{2}\(s_{i}\ou^i(t)\cos\theta_{i+1}(t) -s_{i+1}\ou^{i+1}(t)\cos\theta_{i+1}(t)\) &\text{ if } \cos\theta_{i}(t)=\cos\theta_{i+1}(t) \text{ and } s_{i}\ou^i(t)\neq s_{i+1}\ou^{i+1}(t),\\
0,& \text{ otherwise.}
\end{cases}
$$
The cost functional is calculated by
\begin{equation*}
J[x,u]:=\disp\frac{1}{2}\big\|x(T)\big\|^2=\dfrac{1}{2}\(\(x^{11}(T)\)^2+\(x^{12}(T)\)^2 +\(x^{21}(T)\)^2+\(x^{22}(T)\)^2+\ldots+\(x^{n1}(T)\)^2+\(x^{n2}(T)\)^2\),
\end{equation*}
since $x(t)$ has the same dimension as $x^j_*$ and $c_j$.

\section{Numerical Algorithm}\label{Num}
\setcounter{equation}{0}

This section illustrates the power of the necessary optimality conditions obtained in Theorem~\ref{Thm6.2*a} for solving the controlled robotics model developed in Section~\ref{Prob-Form}. We furnish this by applying the obtained conditions to problems with specified numerical data. Our main achievement here is to design a {\em new algorithm}, which allows us to calculate optimal solutions in robotic models with encompassing different scenarios, where $n=2,\,3,\,4,\,5,\,10,\,100,\ldots$ and beyond.

Based on the sweeping control description of the robotics model in Section~\ref{Prob-Form} and the necessary optimality conditions of Theorem~\ref{Thm6.2*a}, we consider the following optimal control problem for developing our algorithm for in the general robotic setting with an arbitrary number of robots:
$$
\begin{aligned}
&\mbox{ minimize }\;\frac{1}{2}\|x(T)\|^2\;\text{ subject to }\\
-\dot{x}(t)&=\sum_{j=1}^{n-1} \eta_j(t) x_j^*+(g( x_1(t),u_1(t)),\ldots, g( x_n(t),u_n(t)))\\
&\|x_j(t)-x_{j-1}(t)\|\geq 2R
\end{aligned}
$$
with the representation functions $\eta_i$ defined by
$$\eta_j(t)=\begin{cases} 0& \text{if } x_{1,j+1}(t)-x_{1,j}(t)+x_{2,j+1}(t)-x_{2,j}(t)>2R,\\
\disp\frac{1}{2}(s_j u_j(t)\cos\theta_{j+1}(t)-s_{j+1}u_{j+1}(t)\cos\theta_{j+1}(t))& \text{ otherwise}
\end{cases}
$$ 
for $j=1, \ldots, n-1.$ We utilize the GEKKO Dynamic Optimization Suite library from Python, which uses nonlinear solvers such as IPOPT and SNOPT to obtain the numerical solution $(\bar{x}(t),\bar{u}(t))$. 

Our implementation consists of three parts. First, we develop the necessary Python code to efficiently read input files both as a single file and in batch mode to allow for large numbers of simulations.  Second, we implement the code that uses the input file information to formulate the problem as a GEKKO optimization problem that is sent to the selected solver. Finally, the final portion of the code ensures that all the aspects of the simulation procedure, such as input file, output graphics, and solution efficiency data, are recorded and time stamped so that we can analyze large numbers of simulations.

The input files contain crucial parameters such as the initial position of the robots, the bounds on the control function, the robots' radius $R$, and the final time $T$. Leveraging this input, our code initializes the number of robots $(n)$ and executes the GEKKO dynamic optimization algorithm made up of differential equations extracted from the obtained necessary optimality conditions. Successful execution results in a detailed description of the path of the robots, the controls utilized, the performance measure, and the functions $\eta_i$ are presented in Figure~1 below.
\begin{figure}[!h]
 \centering
\includegraphics[width=\textwidth]{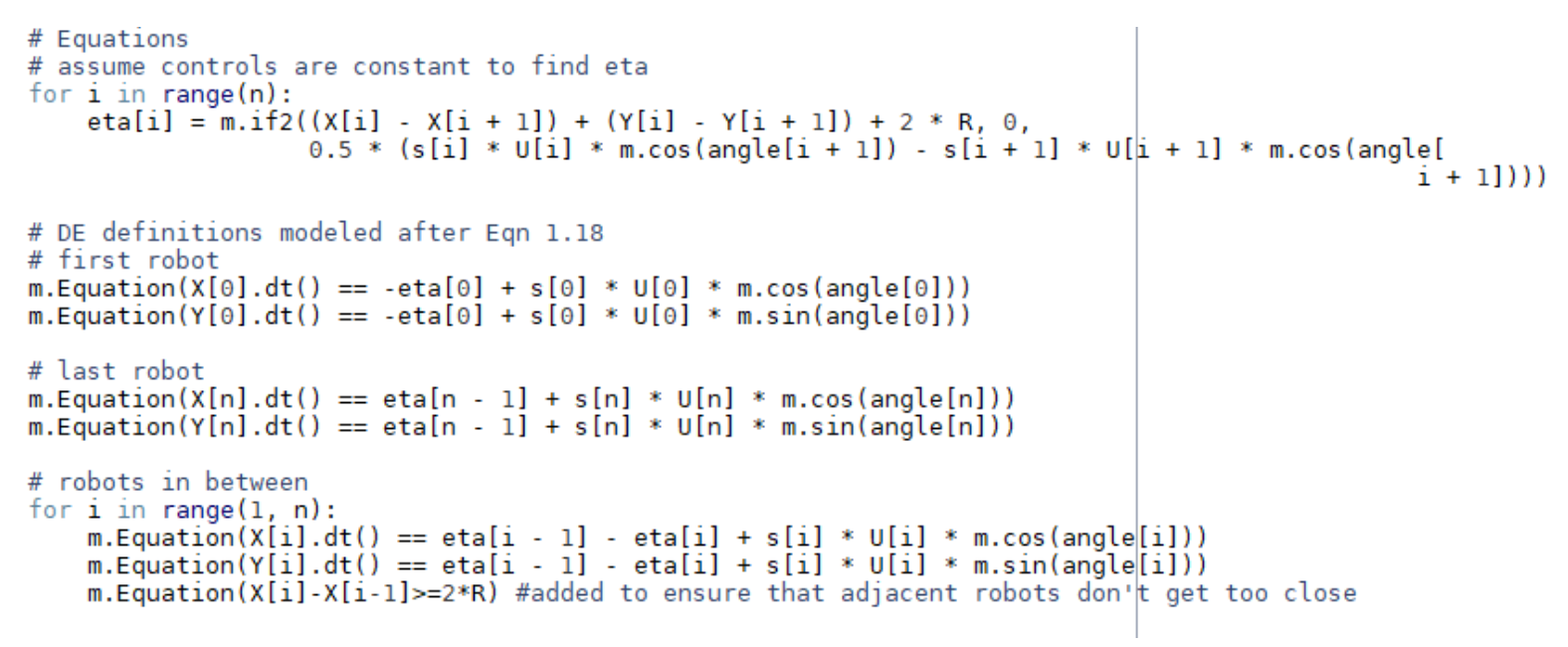}
\caption{The code with the differential equations involved in the
necessary conditions}\label{codepython}
\end{figure}

In our work so far, we have been able to run around $200$ simulations. Although we performed some simulations utilizing the IPOPT solver and others with SNOPT, we did not observe any significant differences in time to a solution or likelihood of convergence based on the solver chosen.  We did find a significant difference in time of convergence based on the number of robots and observed that the time to achieve a solution increased rapidly and proportionally to the number of robots involved.  We illustrate our results here with three examples. The first is a simulation that involves two robots. Solutions to $2-$robot simulations resulted in one of three scenarios. Either one robot moves as close as possible to the origin and then stops to allow the second robot to proceed to the origin as fast as possible, the two move almost in synchrony towards the origin, or finally the robots wait to move and then move rapidly to the origin.  In all the cases, the time to a solution hovered around $10$ seconds and the performance measure was often $0.5R$. This last observation about the performance measure highlights one of the challenges of the problem as stated.  Since the robots are set to reach the origin, but are not removed from the dynamics once they do, in order to minimize the distance of all robots to the origin, the solution must balance the distance of all robots to the origin causing at times for robots to move past the origin. In the case of two robots, this can be achieved by having one robot $0.5R$ past the origin and one the same distance on the opposite side of the origin. We specify the data as follows $x^{11}(0)=2,\;x^{12}(0)=2,\;x^{21}(0)=3,\;x^{22}(0)=3,\;T=2,\;R_1=R_2=1,\;s_1=3,\;s_2=1$, $\theta_i=45^\circ$ and $\|u_i(t)\| \leq 10$, for $i=1,2.$ The controls and the position of the robots are presented below and illustrate the second scenario, where the robots move almost synchronously toward the origin, and where the final position shows both robots at the positions described ($1/2R$).

\begin{center}
\begin{frame}{Two robots.  $R=1$, $T=2$, and $|u(t)|\leq 10$. }
\includegraphics[width=5in]{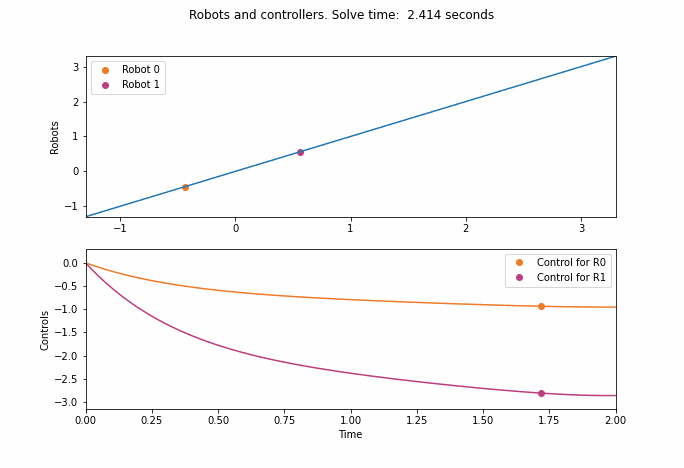}

\includegraphics[width=4in]{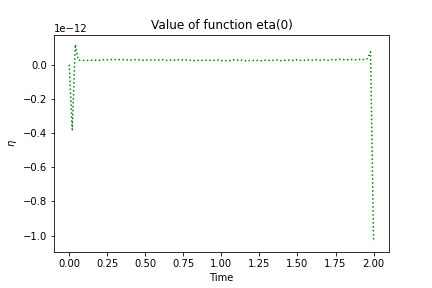}
 \end{frame}
\end{center}

In this second example, we include a simulation with $5$ robots. The data for this case is contained in Table~\ref{tab:values5}, where $(X,Y)$ are the positions of the robots on the plane at the initial time, $s$ indicates the speed for each robot. The total time is $T=8$, five robots have the same radius $R=1$, $\theta_i=45^\circ$ for $i=0,\ldots,4$, and the controls satisfy $\|u_i\|\leq 5$ for $i=0,\ldots,4$.
\begin{table}[!ht]
\centering
\begin{tabular}{|c|c|c|c|c|c|}
\hline
\textbf{Robot} & \textbf{X} & \textbf{Y} & \textbf{s} \\
\hline
0 & 5 & 5 & 2 \\
\hline
1 & 11 & 11 & 2 \\
\hline
2 & 16 & 16 & 1 \\
\hline
3 & 20 & 20 & 3 \\
\hline
4 & 27 & 27 & 3 \\
\hline
\end{tabular}
\caption{Data for the case of 5 robots.}
\label{tab:values5}
\end{table}

The time to solution was $211$ seconds and the performance measure of $32$. On average simulations with $5$ bots took about $200$ seconds.  This example demonstrates the outcome of having two robots wait for the others to get close to the origin before speeding there.  One also can see $\eta$ activated around time $t=6$ when the last two robots are to move, but would collide if allowed to use their natural speeds. 

\begin{center}
\begin{frame}{Five robots simulation. $R=1$, $T=8$, and $|u(t)|\leq 5$. }
 \includegraphics[width=5in]{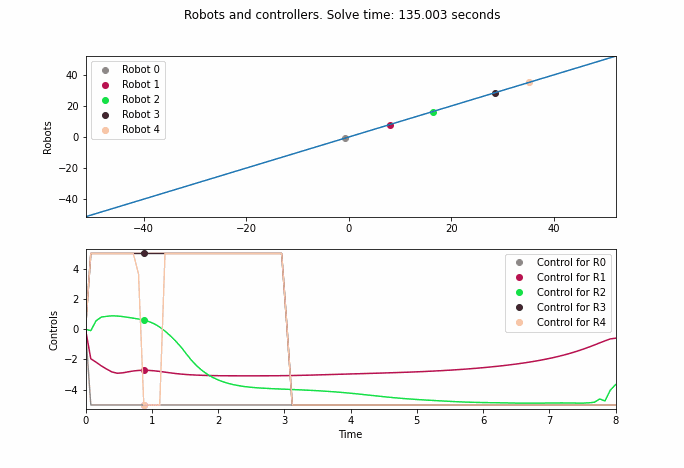}
\end{frame}
\end{center}
\begin{center}
\begin{figure}[htp]
\includegraphics[width=3in]{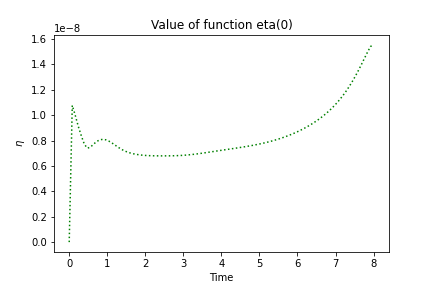}
\includegraphics[width=3in]{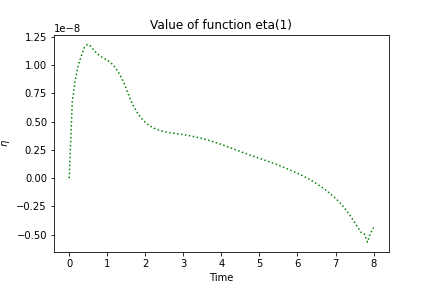}
\includegraphics[width=3in]{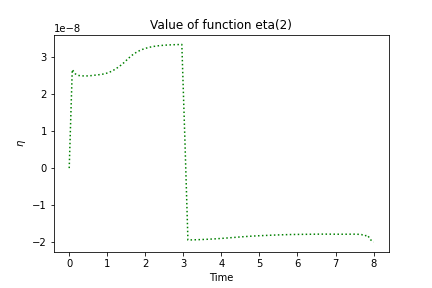}
\includegraphics[width=3in]{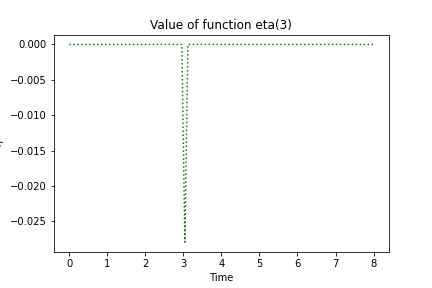}
\end{figure}
\end{center}

The last example is the one with $10$ robots. Table~\ref{tab:values} shows the data for this case with $T=8$, radius for all robots is $R=1$, $\theta_i=225^\circ$, and the controls $\|u_i(t)\|\leq 5$, for $i=0,\ldots,9$.
\begin{table}[!ht]
\centering
\begin{tabular}{|c|c|c|c|}
\hline
\textbf{Robot}& \textbf{X} & \textbf{Y} & \textbf{s} \\
\hline
0 & -10 & -10 & 1 \\
\hline
 1 & -13 & -13 & 1 \\
\hline
 2 & -18 & -18 & 1 \\
\hline
 3 & -22 & -22 & 2 \\
\hline
 4 & -25 & -25 & 2 \\
\hline
 5 & -29 & -29 & 2 \\
\hline
 6 & -32 & -32 & 3 \\
\hline
7 & -35 & -35 & 3 \\
\hline
8 & -38 & -38 & 3 \\
\hline
9 & -40 & -40 & 3 \\
\hline
\end{tabular}
\caption{Data for the case of 10 robots}
\label{tab:values}
\end{table}}
On average these simulations took $1100$ seconds, but the example below took $1751.3$ seconds. In this simulation, we can see that the robots have to adjust their speed almost immediately and then again before the time $T=3$; both cases requiring the activation of $\eta$.  After that they move almost synchronously toward the origin. 
\begin{center}
\begin{frame}{Ten Robots simulation.  $R=1$,$T=8$ and $|u(t)|\leq 5$. }
\includegraphics[width=5in]{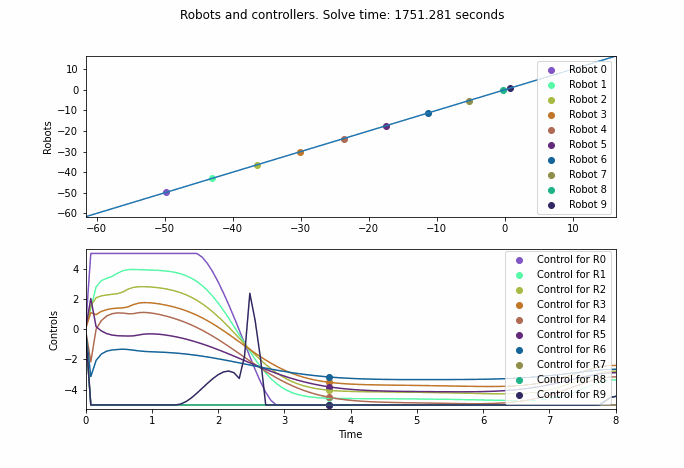}
\end{frame}
\end{center}
\begin{center}
\begin{figure}[htp]
\includegraphics[width=3in]{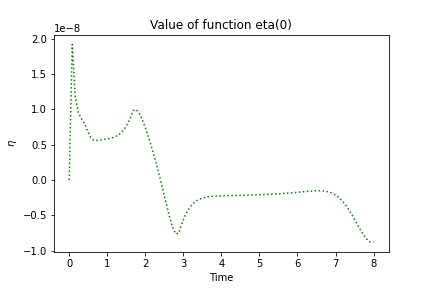}
\includegraphics[width=3in]{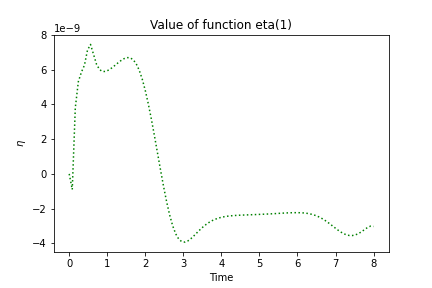}
\includegraphics[width=3in]{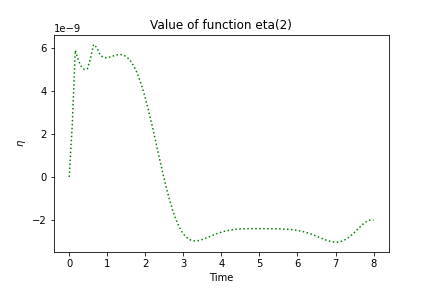}
\includegraphics[width=3in]{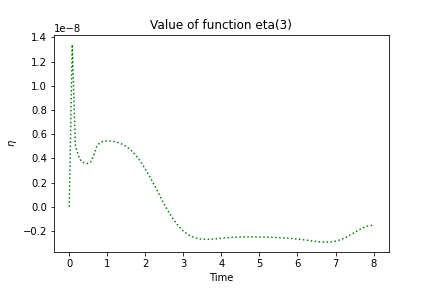}
\includegraphics[width=3in]{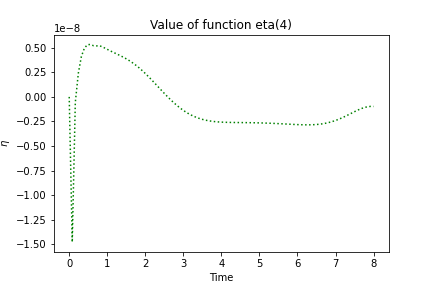}
\includegraphics[width=3in]{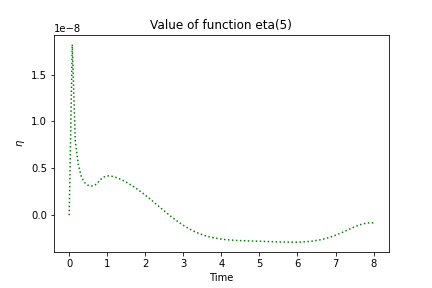}
\includegraphics[width=3in]{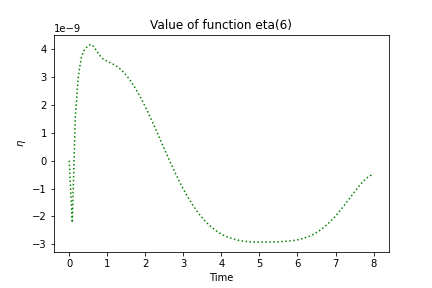}
\includegraphics[width=3in]{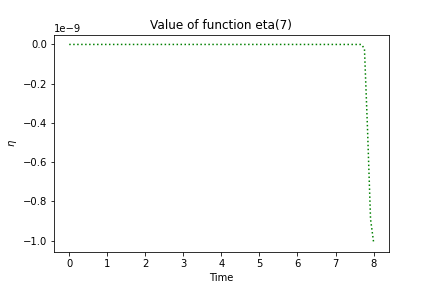}
\end{figure}
\end{center}
\begin{center}
\includegraphics[width=3in]{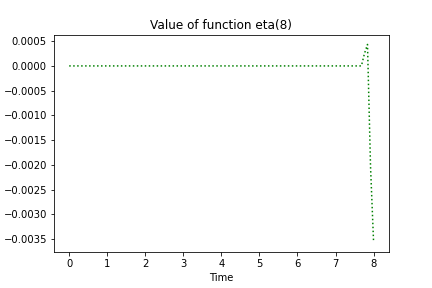}
\end{center}

These simulations are promising, but there are important adjustments that may be considered and that may affect the efficiency of the solvers. For instance, the code could account for success when the robot reaches the origin in order that the final solution is not a minimizing weighted position for all the robots as related to the origin, with some unable to reach the target.  Implementing such adjustments would result in smaller performance measures, and it is possible that they could lead to faster convergence times. Overall, our method of solving reaches good visualizations; however, in its current state, the time to a solution is not reasonable for very large numbers of robots.  One possible explanation for this limitation is that the time allocated for the robots to reach the target may not be optimal. Therefore, our future work will be devoted to  minimizing the time variable required for the robots to reach the target in the algorithm, which takes into account the theoretical results developed in \cite{cmnn23b} with applications to other practical sweeping control models, e.g., those considered in \cite{mnn23a}.

 {\bf Acknowledgements.} The authors are grateful to Andrew Regan for many useful discussions related to the coding part of the algorithm.

\end{document}